\documentclass[11pt]{amsart}  
\usepackage[left=1.3in,right=1.3in,bottom=1.3in]{geometry}  

\usepackage{amssymb,amsmath,amscd}
\usepackage{amsthm}

\usepackage{mathrsfs}
\usepackage{latexsym}
\usepackage{enumerate}
\usepackage{graphicx} 
\usepackage{tikz-cd}
\usepackage{galois}
\usepackage{stmaryrd}
\usepackage{algorithm,algorithmic}

\theoremstyle{plain}
\newtheorem{theorem}[]{Theorem}

\newtheorem{corollary}[theorem]{Corollary}
\newtheorem{definition}[]{Definition}

\theoremstyle{definition}
\newtheorem{example}{Example}
\newtheorem{remark}{Remark}

\newcommand{\va}{\textsc{val}}
\newcommand{\ca}{\textsc{cap}}
\newcommand{\topL}{\overline{\ell}}
\newcommand{\botL}{\underline{\ell}}

\newcommand{\style}[1]{{\bf #1}}

\begin{document}

\title{Lattice-Valued Bottleneck Duality}
\author{Robert Ghrist}
\author{Julian Gould}
\author{Miguel Lopez}
\address{Department of Mathematics, Department of Electrical \& Systems Engineering, and Program in Applied Mathematics \& Computational Science,  University of Pennsylvania}


\begin{abstract}
This note reformulates certain classical combinatorial duality theorems in the context of order lattices. For source-target networks, we generalize bottleneck path-cut and flow-cut duality results to edges with capacities in a distributive lattice. For posets, we generalize a bottleneck version of Dilworth's theorem, again weighted in a distributive lattice. These results are applicable to a wide array of non-numerical network flow problems, as shown. All results, proofs, and applications were created in collaboration with AI language models. An appendix documents their role and impact.
\end{abstract}

\maketitle

\section{Introduction}

The Max-Flow Min-Cut (MFMC) theorem is a cornerstone of network flow theory and establishes a duality between the maximum flow value that can be pushed through a network and the minimum capacity of a cut separating source from target. It is a simple consequence of LP-duality and, as with linear programming, is of profound importance. Bottleneck duality, though more restrictive and less prominent, is nevertheless an important example of duality --- path-cut instead of flow-cut.

The modest goal of this paper is to extend bottleneck duality to systems valued in an order lattice. The perspectives and language of lattice theory will be used throughout: see \cite{davey2002introduction,gratzer2011lattice} for relevant background. 

\subsection{Classical flow-cut duality theorems}

\begin{theorem}[MFMC]
\label{thm:MFMC}
Let $G = (V, E)$ be a finite, connected directed graph without self-loops or parallel edges. Let $s, t \in V$ be distinct source and sink vertices, respectively, and let $c: E \to \mathbb{R}^+$ a capacity function on edges. Then:
\[
    \max_{\phi\in\Phi}
    \va(\phi)
    = 
    \min_{C\in\mathcal{C}}
    \ca(C) 
\]
where:
\begin{itemize}
\item $\Phi$ is the (non-empty) set of feasible flows from $s$ to $t$,
\item $\mathcal{C}$ is the (non-empty) set of all cuts separating $s$ from $t$,
\item $\va(\phi) = \sum_{v:(s,v) \in E} \phi(s,v)$ denotes the value of a flow $\phi$,
\item $\ca(C) = \sum_{e \in C} c(e)$ denotes the capacity of an $s$-$t$ cut $C$.
\end{itemize}
Here, a feasible flow $\phi: E \to \mathbb{R}^+$ satisfies:
\begin{enumerate}
\item Capacity constraints: $\forall e \in E, \; 0 \leq \phi(e) \leq c(e)$
\item Flow conservation: $\forall v \in V \setminus \{s,t\}, \; \sum_{u:(u,v) \in E} \phi(u,v) = \sum_{w:(v,w) \in E} \phi(v,w)$
\end{enumerate}
\end{theorem}

Shifting focus from the total flow through the network to the throughput of individual paths converts flow values to path bottlenecks --- the edge along a path with minimum capacity. The corresponding result is:

\begin{theorem}[Bottleneck Duality]
\label{thm:BD}
Let $G = (V, E)$ be a finite, connected directed graph without self-loops or parallel edges. Let $s, t \in V$ be distinct source and sink vertices, respectively, and let $c: E \to \mathbb{R}^+$ be a capacity function on edges. Then:
\[
\max_{P \in \mathcal{P}} \min_{e \in P} c(e) = \min_{C \in \mathcal{C}} \max_{e \in C} c(e)
\]
where $\mathcal{P}$ is the (non-empty) set of all paths from $s$ to $t$, and $\mathcal{C}$ is the (non-empty) set of all $s$-$t$ cuts.
\end{theorem}
We will generalize this result to capacities taking values in distributive lattices.

The classic MFMC theorem dates back to the seminal work in the late 1950s \cite{ford1956maximal,elias1956note,dantzig1956max}. The origins of the bottleneck flow problem can be traced back to the early 1960s \cite{Pollack1960,Hu1961}. 

Manifestations of duality in combinatorial optimization appear beyond the classical and bottleneck max-flow min-cut theorems. In graph theory, Menger's Theorem \cite{menger1927allgemeinen} establishes a duality between the maximum number of vertex-disjoint paths and the minimum vertex cut in a graph, presaging the more general flow-cut duality concepts. K{\"o}nig's Theorem \cite{konig1931graphen} reveals a duality between maximum matching and minimum vertex cover in bipartite graphs, which can be viewed as a specialized instance of flow-cut duality. Dilworth's Theorem \cite{dilworth1950decomposition} extends these ideas to partially ordered sets, establishing a duality between the maximum size of an antichain and the minimum number of chains needed to cover the set.

The Multicommodity Flow-Cut Theorem of Leighton-Rao \cite{leighton1999multicommodity} generalizes max-flow min-cut to scenarios involving multiple commodities, linking the maximum concurrent flow to the sparsest cut in the network. The Matroid Intersection Theorem \cite{edmonds1970submodular} further abstracts these concepts, establishing a duality between the maximum size of an intersection and a certain type of partition in the context of matroids. In the domain of combinatorial game theory, the Gale-Berlekamp Switching Game \cite{berlekamp1965design} presents a surprising connection to flow-cut duality. The Lovász-Plummer Matching Forest Theorem \cite{lovasz2009matching} and Nash-Williams Arborescence Theorem \cite{nash1961edge} extend flow-cut duality to more complex graph structures, dealing with matching forests and edge-disjoint spanning trees respectively. 

\subsection{Background}

A \style{lattice} is a partially ordered set $(L, \leq)$ in which every pair of elements $a, b \in L$ has a unique supremum (least upper bound) called the \style{join}, denoted $a \vee b$, and a unique infimum (greatest lower bound) called the \style{meet}, denoted $a \wedge b$. The operations $\vee$ and $\wedge$ are required to satisfy:
\begin{enumerate}
\item Commutativity: $a \vee b = b \vee a$ and $a \wedge b = b \wedge a$;
\item Associativity: $a \vee (b \vee c) = (a \vee b) \vee c$ and $a \wedge (b \wedge c) = (a \wedge b) \wedge c$;
\item Absorption: $a \vee (a \wedge b) = a$ and $a \wedge (a \vee b) = a$;
\item Idempotence: $a \vee a = a$ and $a \wedge a = a$.
\end{enumerate}

The ordered reals are a lattice with max and min as join and meet. Other common examples of lattices include: the power set of a set ordered by inclusion, with union as join and intersection as meet; Boolean algebras ordered by implication, with disjunction as join and conjunction as meet; and the set of partitions of a set ordered by refinement, with the coarsest common refinement as join and the finest common coarsening as meet. A wealth of useful lattices can be found in the literature on Formal Concept Analysis \cite{ganter1999formal,wille1982restructuring}.

A lattice is \style{complete} if every subset (including the empty set and the entire lattice) has both a join and a meet. 
A \style{modular} lattice satisfies the modular law: for all $a, b, c \in L$, if $a \leq c$, then $a \vee (b \wedge c) = (a \vee b) \wedge c$.
A lattice is \style{distributive} if binary joins distribute over binary meets and vice versa:
\begin{equation}
\label{eq:distributive}
    a\vee(b\wedge c) = (a\vee b)\wedge(a\vee c)
    \quad \text{and} \quad
    a\wedge(b\vee c) = (a\wedge b)\vee(a\wedge c) .
\end{equation}
All distributive lattices are modular. 

We will work with the following class of networks which have a specified source and target.
\begin{definition}[Flow Network]
\label{def:flownet}
A {\bf flow network} is a tuple $(G, s, t)$ where:
\begin{enumerate}
\item $G = (V, E)$ is a finite, directed graph without self-loops or directed cycles; 
\item $s, t \in V$ are distinct vertices called the source and sink, respectively;
\item All edges incident to $s$ are outgoing, and all edges incident to $t$ are incoming.
\item For all $v \in V \setminus \{s, t\}$, there exist directed paths $s \rightsquigarrow v$ and $v \rightsquigarrow t$ in $G$.
\end{enumerate}
\end{definition}
Some of our results extend to more general networks, such as those with dead-ends: see Remarks \ref{rem:deadend1} and \ref{rem:deadend2}. 

A \style{cut} $C = (S, T)$ is a partition of $V$ such that $s \in S$ and $t \in T$. By way of abuse of terminology, one says that $e\in C$, for an edge $e=(u,v)$ if $u\in S$ and $v\in T$. A cut $C$ is \style{minimal} if there are no cuts $C'$ such that $\{e \in E \: : \: e \in C' \} \subsetneq \{e \in E \: : \: e \in C \}$.

\section{Bottleneck duality}

All results in this note stem from bottleneck duality. 

\subsection{Result}

The following is a lattice-valued generalization of Theorem \ref{thm:BD}.

\begin{theorem}[Lattice Bottleneck Duality]
\label{thm:LVBD}
For $G$ a flow network with $c: E \to L$ a capacity function taking values in a distributive lattice $L$:
\begin{equation}
\label{eq:LVBD}
        \bigvee_{P\in\mathcal{P}} \bigwedge_{e \in P} c(e) 
        = 
        \bigwedge_{C\in\mathcal{C}} \bigvee_{e \in C} c(e) ,
\end{equation}
where $\mathcal{P}$ is the set of all paths from $s$ to $t$, and $\mathcal{C}$ is the set of all cuts separating $s$ from $t$. 
\end{theorem}
{\em Proof:} 
Let $\alpha = \bigvee_{P\in\mathcal{P}} \bigwedge_{e\in P} c(e)$ and $\beta = \bigwedge_{C\in\mathcal{C}} \bigvee_{e \in C} c(e)$. 

We begin with weak duality: $\alpha\leq\beta$.
Consider any path $P \in \mathcal{P}$ and any cut $C \in \mathcal{C}$. Since $P$ goes from $s$ to $t$ and $C$ separates $s$ from $t$, there must be at least one edge $e' \in P \cap C$
. 
For this edge, we have:
\[
\bigwedge_{e\in P} c(e) \leq c(e') \leq \bigvee_{e\in C} c(e).
\]
The first inequality holds because $e' \in P$, and the meet of a set is less than or equal to any element of the set. The second inequality holds because $e' \in C$, and any element of a set is less than or equal to the join of the set.

Since this holds for any $P \in \mathcal{P}$ and any $C \in \mathcal{C}$, we can take the join over all $P$ on the left side and the meet over all $C$ on the right side:
\[
\alpha = \bigvee_{P\in\mathcal{P}} \bigwedge_{e\in P} c(e) \leq \bigwedge_{C\in\mathcal{C}} \bigvee_{e \in C} c(e) = \beta.
\]

Strong duality, $\beta \leq \alpha$, is more delicate. 
In a distributive lattice, finite meets distribute over finite joins and vice versa. For any array of elements $\{a_{i,j}\} \subset L$ over finite index sets $I$ and $J$,
\begin{equation}
\label{eq:findist}
\bigwedge_{i \in I} \bigvee_{j \in J} a_{i,j} = \bigvee_{j \in J} \bigwedge_{i \in I} a_{i,j} .
\end{equation}
We will apply this to our context with paths and cuts as (finite) index sets. Define $a_{C,P}$ for $C \in \mathcal{C}$ and $P \in \mathcal{P}$ as:
\[
    a_{C,P} = \bigvee_{e \in C \cap P} c(e) .
\]
Note that $a_{C,P}$ is well-defined because $G$ is finite, thus $C \cap P$ is finite, and the join is over a finite (and non-empty) set.

\begin{equation}\label{eq:LVBNSD}
    \beta = \bigwedge_{C \in \mathcal{C}} \bigvee_{e \in C} c(e) 
    \leq \bigwedge_{C \in \mathcal{C}} \bigvee_{P \in \mathcal{P}} a_{C,P} = \bigvee_{P \in \mathcal{P}} \bigwedge_{C \in \mathcal{C}} a_{C,P} 
    \leq 
    \bigvee_{P \in \mathcal{P}} \bigwedge_{e \in P} c(e) 
    = \alpha
\end{equation}

The first equality is the definition of $\beta$. The second inequality holds because every edge in a cut is on some path via property (4) in Definition \ref{def:flownet}. 
The third equality is justified by finite meets distributing over finite joins in a distributive lattice. The fourth inequality holds because for every $P \in \mathcal{P}$ and $e \in P$, there is a cut $C \in \mathcal{C}$ such that $P \cap C = \{e\}$. Hence $\bigwedge_{C \in \mathcal{C}} a_{C,P} \leq \bigwedge_{e \in P} c(e)$ and so the fourth inequality follows. The final equality is the definition of $\alpha$.

This proves strong duality. Combining weak and strong duality demonstrates that $\alpha = \beta$.
\qed

\begin{remark}\label{rem:deadend1}
    If the network has dead ends, Theorem \ref{thm:LVBD} still holds assuming the lattice $L$ has a minimal element zero. Some care must be taken in the proof to restrict to minimal cuts when proving the second inequality in Equation \eqref{eq:LVBNSD}. Finally, in the case there are no paths from $s$ to $t$, the conclusion of Bottleneck Duality is merely that $0 = 0$. 
\end{remark}

\subsection{Examples}

The classical $\mathbb{R}^+$-valued bottleneck duality follows from using the lattice $L=(\mathbb{R}^+,\leq)$ with $\wedge$ as min and $\vee$ as max. Other examples demonstrate the limits of our results.

\begin{figure}[ht]
\centering
    \includegraphics[width=5in]{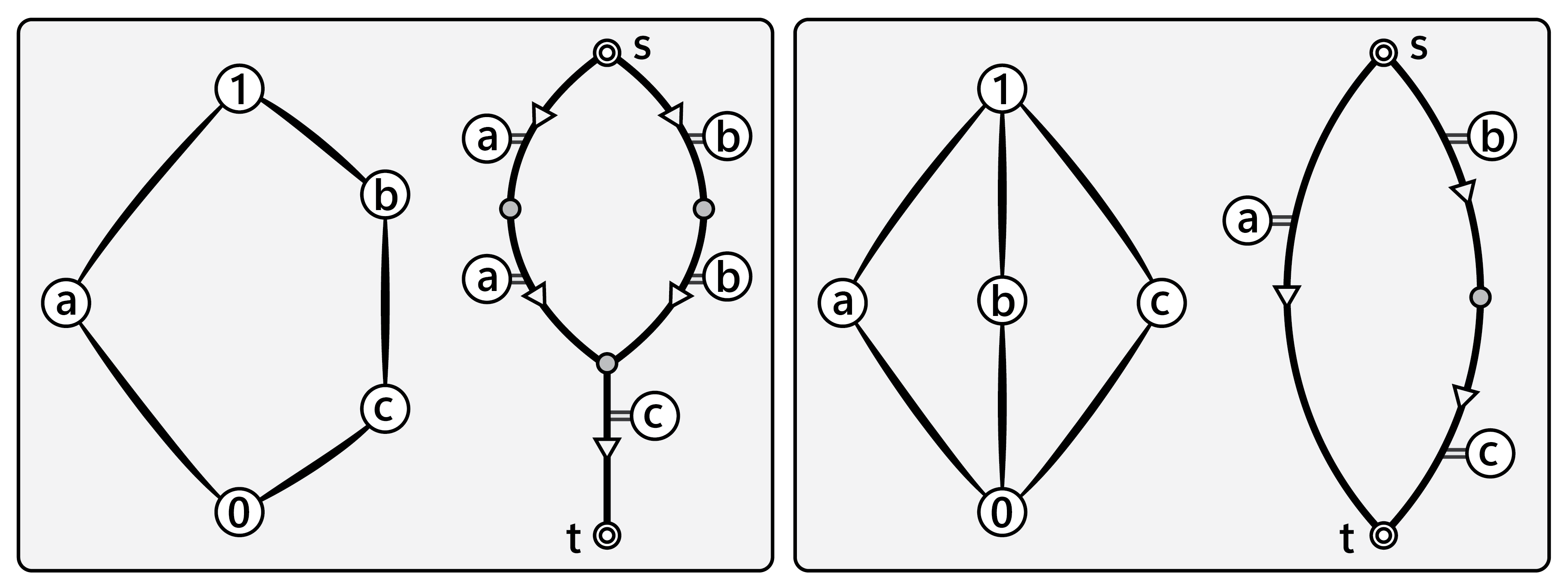}
\caption{Bottleneck duality fails for [left] the pentagon lattice on a simple flow network with five edges, and [right] the diamond lattice on a simple flow network with three edges.}
\label{fig:1}
\end{figure}

\begin{example}
\label{ex:pentagon}
 Let $L$ be the {\em pentagon lattice} with elements $\{0, a, b, c, 1\}$, where $0 < c < b < 1$, $0 < a < 1$, and $a$ is incomparable to both $b$ and $c$: see Figure \ref{fig:1} [left]. The simple flow network there illustrated has 
 two paths from $s$ to $t$ with throughputs $a \wedge b = 0$, and $c \wedge b = c$. There are six cuts separating $s$ from $t$ with capacities $a \vee c = 1$, and $b$. Bottleneck duality would assert that $0\vee c = 1\wedge b$: false.
 The lattice $L$ is not distributive and bottleneck duality fails.
\end{example}

\begin{example}
\label{ex:diamond}
One could hope that modularity would suffice, but this is not the case. Figure \ref{fig:1}[right] illustrates the {\em diamond lattice} with elements $\{0, a, b, c, 1\}$, with $a, b, c$ incomparable. A simple flow network with capacities as labeled has two paths (throughputs $a$ and $b\wedge c=0$) and two cuts (each with capacity $1$). Bottleneck duality clearly fails here as well.
\end{example}

\begin{corollary}
\label{cor:iff}
Equation (\ref{eq:LVBD}) holds for all flow networks with capacities valued in a lattice $L$ if and only if $L$ is distributive.
\end{corollary}
{\em Proof:} 
Theorem \ref{thm:LVBD} is one direction; Examples \ref{ex:pentagon} and \ref{ex:diamond} show the reverse direction, since a lattice is distributive if and only if it does not have a sublattice isomorphic to the pentagon or diamond lattices \cite{birkhoff1967lattice}. 
\qed

\begin{figure}[ht]
\centering
    \includegraphics[width=4.5in]{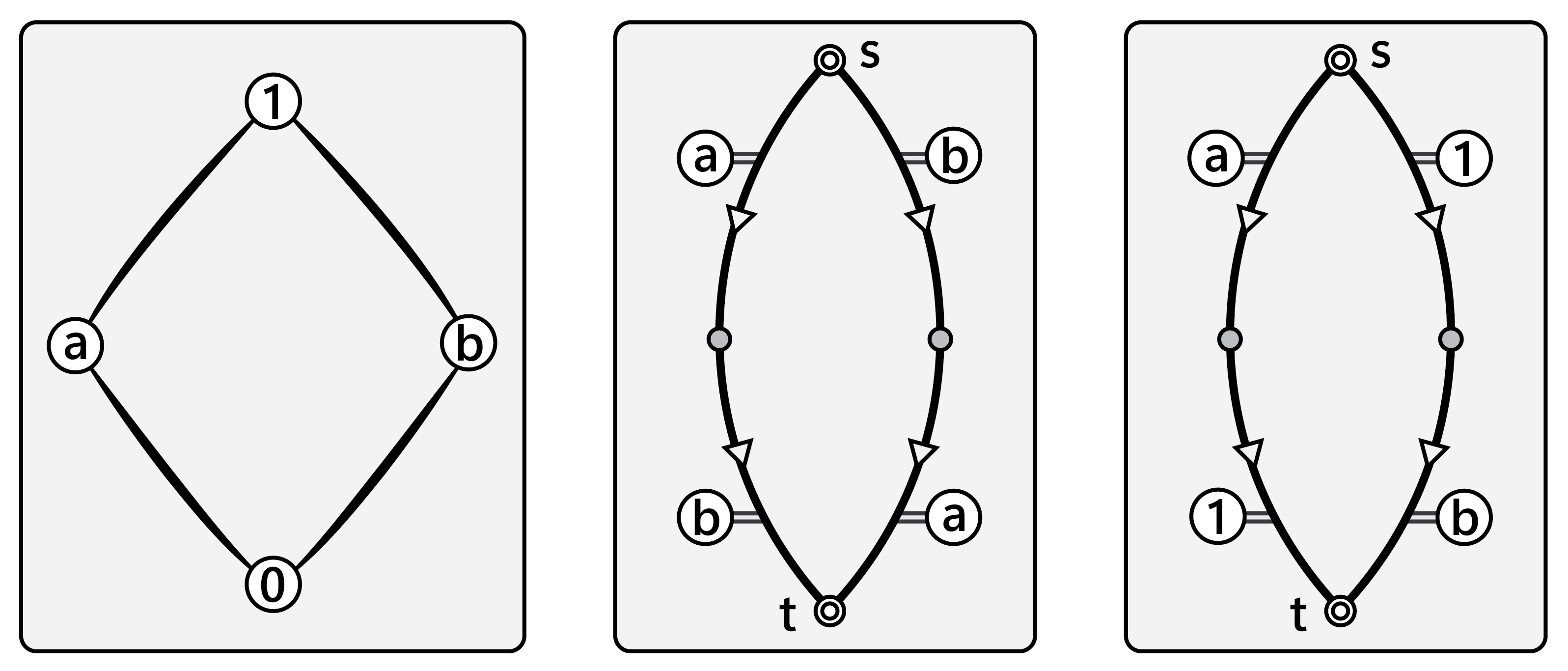}
\caption{A lattice with four elements $L=\{0, a, b, 1\}$, isomorphic to a powerset lattice on $\{a, b\}$, is distributive [left]. A flow network may have optimal paths but no optimal cut [center] and/or optimal cuts but no optimal paths [right].}
\label{fig:2}
\end{figure}

\begin{example}
\label{ex:optimal}
Theorem \ref{thm:LVBD} is a statement on the duality of values, not a statement about existence. Optimal paths and optimal cuts may or may not exist. Figure \ref{fig:2} illustrates two examples of capacities on a flow network taking values in a simple distributive lattice [left], showing that optimal paths and cuts have independent existence and non-existence. In the first flow network [center], all paths have (optimal) throughput $0$, but all cuts have value strictly greater than $0$, though their meet is $0$. In the second example [right], the two paths have throughputs $a$ and $b$, whose join is $1$; yet all cuts have (optimal) cut value $1$. 
The case of a flow network with neither optimal path nor optimal cut is left to the reader. 
\end{example}

\subsection{Applications}

Our goal is to show that the use of lattice coefficients greatly expands the applicability of bottleneck duality to use-cases beyond scalar values. 

\subsubsection{Supply Chain Resource Allocation}
\label{app:supply-chain}

We begin with a simple use of a powerset lattice in the domain of supply chain management, particularly in resource allocation. By modeling a supply chain as a flow network $G = (V, E)$, we can represent various entities such as suppliers,  distributors, and retailers as vertices, with the edges denoting the pathways through which resources flow. 

A typical system might have multiple sources (representing different material source supply) and multiple targets (each of whom might receive a certain subset of resources). It is natural to add an artificial ``master'' source node $s$ and target node $t$ to set up a flow network.

Let $R$ be the finite set of all distinct resources. The lattice $L = 2^R$ then represents all possible combinations of these resources. This lattice structure $(L, \subseteq, \vee, \wedge)$ is naturally ordered by the subset relation, with union and intersection serving as the join and meet operations, respectively.

We assign a capacity function $c: E \to L$ to each edge in the network, where $c(e)$ specifies the subset of resources available for allocation along that edge. This formulation allows us to model the flow of multiple resources simultaneously through the supply chain.

In this context, paths from source to sink represent specific routes through which resources flow, while cuts partition the network into two disjoint subsets, separating the source from the target. The bottleneck value of a path $P$ is given by $c(P) = \bigwedge_{e \in P} c(e) = \bigcap_{e \in P} c(e)$, representing the common resources available across all edges in the path. Conversely, the capacity of a cut $C$ is defined as $c(C) = \bigvee_{e \in C} c(e) = \bigcup_{e \in C} c(e)$, indicating the aggregate resources available across all edges crossing the cut.

Applying the Lattice Bottleneck Duality theorem to this supply chain model yields:
\[
\bigcup_{P \in \mathcal{P}} \left(\bigcap_{e \in P} c(e)\right) = \bigcap_{C \in \mathcal{C}} \left(\bigcup_{e \in C} c(e)\right)
\]
This duality has the following straightforward interpretation:
\begin{quote}
{\em The most comprehensive set of resources consistently available along at least one complete supply path is identical to the minimal set of resources that must be present across every possible critical partition of the supply chain.}
\end{quote}

The left-hand side of the equation represents the aggregate of limiting resource sets across all possible supply paths, identifying the broadest possible set of resources that can be guaranteed to flow through the system via at least one path. The right side represents the intersection of the aggregate resource sets across all possible critical partitions of the supply chain. This identifies the core set of resources that must be available at every potential bottleneck in the system.

\subsubsection{Packaging and Safety}
\label{app:safety}

One can capture more supply chain complexities using the following theorem of Birkhoff.
\begin{theorem}[Birkhoff \cite{birkhoff1967lattice} Ch. 9 Thm. 6]
\label{thm:rep-thm}
        A lattice $L$ is distributive if and only if it is isomorphic to a ring of sets, i.e. a family of sets closed under finite unions and intersections, partially ordered by containment.   
\end{theorem}

Suppose certain resources can only be transported safely with the correct packaging materials. For example, food resources resources might need to be packaged with refrigeration devices, and hazardous materials might need certain insulators. Let $R$ and $M$ denote our sets of resources and packaging materials respectively. To each $r \in R$, let $M_r \subseteq M$ denote the packaging materials required to transport $r$. Taking the closure of $\big \{ \{r \} \cup M_r \: : \: r \in R \big \}$ in $2^{R \sqcup M}$ with respect to finite unions and intersections produces a ring of sets, and hence a distributive lattice $L$ via Theorem \ref{thm:rep-thm}.

With such, the capacity $c(P)$ of a path $P$ from $s$ to $t$ captures precisely which resources can be transported safely with proper packaging materials. Bottlenecks respect the packaging requirements; hence, restricting to this sublattice extends the scope of the model.

\subsubsection{Regulatory compliance}
\label{app:regulatory}

The following uses a lattice of ``fuzzy'' intervals, as per the fuzzy MFMC theorem in \cite{diamond2001fuzzy}. Consider a corporation navigating a complex regulatory landscape to launch a new product. The regulatory process can be modeled as a flow network where edges represent different aspects of compliance, each with a range of achievable compliance levels (topologized as a scalar for simplicity).

Let $G = (V, E)$ be a flow network representing potential regulatory compliance processes, with source $s$ representing the initial product concept and sink $t$ representing full regulatory approval. Multiple varied interconnected intermediate steps lie between source and target. Define a capacity function $c: E \to L$, where $L$ is the lattice of closed intervals $[a, b] \subset [0,1]$ with $a \leq b$, ordered as follows:

For $A = [a_1, a_2]$ and $B = [b_1, b_2]$,
\[
A \leq B \iff a_1 \leq b_1 \text{ and } a_2 \leq b_2
\]

The lattice operations are defined as:
\begin{align*}
A \wedge B &= [\min\{a_1, b_1\}, \min\{a_2, b_2\}] \\
A \vee B &= [\max\{a_1, b_1\}, \max\{a_2, b_2\}]
\end{align*}
With these operations, $(L, \leq, \vee, \wedge)$ forms a distributive lattice.
For each edge $e \in E$, $c(e) = [l_e, u_e]$ represents the range of compliance levels the company can achieve for that step. The lower bound $l_e$ represents the minimum level of compliance the company is certain it can achieve, while the upper bound $u_e$ represents the maximum level it might achieve with additional effort or resources.

Applying Theorem \ref{thm:LVBD} to this regulatory compliance model yields:
\[
\bigvee_{P \in \mathcal{P}} \left(\bigwedge_{e \in P} c(e)\right) = \bigwedge_{C \in \mathcal{C}} \left(\bigvee_{e \in C} c(e)\right)
\]
where $\mathcal{P}$ is the set of all compliance paths from initial concept to full approval, and $\mathcal{C}$ is the set of all cuts separating these states.

{\em Path interpretations:} For a compliance path $P$, $\bigwedge_{e \in P} c(e)$ represents the weakest compliance interval in that compliance strategy. The left-hand side of the equation thus represents the strongest among these weakest links across all compliance strategies.

{\em Cut interpretations:} A cut $C$ is a regulatory obstruction which must be crossed. The cut value, $\bigvee_{e \in C} c(e)$, represents the best possible compliance interval achievable through any route crossing that regulatory boundary. The right-hand side of the equation represents the weakest among these best possible levels across all regulatory boundaries.

Theorem \ref{thm:LVBD} thus implies:
\begin{quote}
{\em
    The strongest among the weakest regulatory intervals across all compliance strategies is identical to the weakest among the best possible compliance intervals across all regulatory boundaries.
}
\end{quote}

\section{Bottleneck max-flow-min-cut}

While we do not convert the classical max-flow-min-cut theorem to a lattice-valued versions, we can adapt lattice-valued bottleneck duality from path-cut to flow-cut duality by redefining flow conservation using the join operator.

\subsection{Result}
\label{sec:LVMFMCresult}

Let $G$ be a flow network equipped with a capacity function $c:E \rightarrow L$. Let us define an $L$-valued \style{flow} to be a function $\phi : E \to L$ such that for each $e \in E$, $\phi(e) \leq c(e)$, and for each vertex $v \in V \setminus \{s, t\}$,
\[
\bigvee_{(u, v) \in E} \phi(u, v) = \bigvee_{(v, w) \in E} \phi(v, w).
\]
The value of a flow, $\phi$, and the capacity of a cut, $C$, are defined as 
\[
\va(\phi) 
    = 
    \bigvee_{(s, v) \in E} \phi(s, v)
    \quad
    :
    \quad
\ca(C) 
    = \bigvee_{e\in C} c(e).
\]

Our use of $\vee$ in the conservation and flow value equations permits translation of bottleneck duality to flow-cut duality.

\begin{theorem}[Lattice-Valued Bottleneck Max-Flow Min-Cut]
\label{thm:LVBMFMC}
Let $G$ be a flow network with capacity $c: E \to L$ from edges to a distributive lattice $L$.
Then the maximal flow value from $s$ to $t$ equals the minimal cut value:
\begin{equation}
\label{eq:LMFMC}    
    \bigvee_{\phi \in \Phi} \va(\phi) = \bigwedge_{C \in \mathcal{C}} \ca(C),
\end{equation}
where $\Phi$ is the set of all $L$-valued flows on $G$, $\mathcal{C}$ is the set of all cuts separating $s$ from $t$.
\end{theorem}

{\em Proof:} For any path $P\in\mathcal{P}$, define $c(P) = \bigwedge_{e \in P} c(e)$. We will show that
\[
\bigvee_{\phi \in \Phi} \va(\phi) = \bigvee_{P \in \mathcal{P}} c(P) = \bigwedge_{C \in \mathcal{C}} \bigvee_{e \in C} c(e) = \bigwedge_{C \in \mathcal{C}} \ca(C).
\]
For each path $P\in\mathcal{P}$, define the flow $\phi_P: E \to L$ as:
\[
\phi_P(e) = \begin{cases} 
c(P) & \text{if } e \in P \\
\botL & \text{otherwise},
\end{cases}
\]
where $\botL = \bigwedge_{e \in E} c(e)$ is the minimal element generated by $c$ in $L$. 
We now verify that $\phi_P$ is indeed a flow:

(1) Capacity constraints: For all $e \in E$, $\phi_P(e) \leq c(e)$.
    If $e \in P$, then $\phi_P(e) = c(P) = \bigwedge_{e' \in P} c(e') \leq c(e)$. If $e \notin P$, then $\phi_P(e) = \botL \leq c(e)$. In both cases, the capacity constraint is satisfied.

(2) Flow conservation: For any $v \in V \setminus \{s,t\}$,
    $\bigvee_{e \in \text{in}(v)} \phi_P(e) = \bigvee_{e \in \text{out}(v)} \phi_P(e)$.
    This holds because for $v \in P$, both sides equal $c(P)$, and for $v \notin P$, both sides equal $\botL$. This confirms flow conservation, and that $\phi_P$ is a valid flow. 

The value of this flow is $\va(\phi_P) = \bigvee_{e \in \text{out}(s)} \phi_P(e) = c(P)$. From this construction, we conclude 
$\bigvee_{\phi \in \Phi} \va(\phi) = \bigvee_{P \in \mathcal{P}} c(P)$, which equals $\bigwedge_{C \in \mathcal{C}} \bigvee_{e \in C} c(e)$ via bottleneck  duality.
\qed

\begin{remark}
\label{rem:deadend2}
    If the network has dead-ends, Theorem \ref{thm:LVBMFMC} still holds assuming the lattice $L$ has a minimal element $0$, which can be used in place of $\botL$ in the definition of $\phi_P$. Conservation still holds at dead ends since the bottom element $0$ is the empty join.
\end{remark}

\subsection{Applications}

Though flow-cut and path-cut bottleneck dualities are basically the same phenomenon, there are a few instances where it is perhaps more natural to work in the context of lattice-valued flows. 

\subsubsection{Secure Information Flows}
\label{app:secureinfo}

Consider an intelligence network where information security encompasses multiple dimensions, such as confidentiality and reliability.  We model these security requirements using a lattice $L = I_C \times I_R$, where $I_C$ and $I_R$ are totally ordered intervals representing confidentiality and reliability levels, respectively. These intervals can be either discrete (e.g., typical clearance levels) or continuous (e.g., quantitative measures of information sensitivity or source reliability). Similar lattice models of secure information  passage were considered in \cite{Denning76} without any assumptions of distributivity.

Let $G = (V, E)$ be a flow network representing the intelligence network, with source $s$ and target $t$. Define a capacity function $c: E \to L$ where $c(e) = (c_C(e), c_R(e))$ denotes the maximum confidentiality and reliability levels that can be securely transmitted through edge $e$.
The partial order $\leq$ is the product of total orders on $I_C$ and $I_R$, and the lattice operations are component-wise min ($\wedge$) and max ($\vee$).

An $L$-flow $\phi: E \to L$ assigns to each transmission edge $e$ a security level $\phi(e) = (\phi_C(e), \phi_R(e))$ such that $\phi(e) \leq c(e)$ for all $e \in E$. Flow conservation ensures:
\[
\bigvee_{(u,v) \in E} \phi(u,v) = \bigvee_{(v,w) \in E} \phi(v,w) \quad \forall v \in V \setminus \{s,t\}.
\]
{\em Flow interpretations:} The flow structure is particularly appropriate for such an intelligence network, as each vertex represents a ``router'' that must be capable of both receiving and transmitting information at the same maximal security levels. This flow condition, as opposed to a path condition, ensures that each node in the network can handle the incoming information and redistribute it without compromising security levels.

The flow value represents the minimal levels at which intelligence can be routed through the entire system without violating any security constraints. It captures the highest combination of confidentiality and reliability that can be consistently maintained across the network.

{\em Cut intepretations:} In this context, cuts represent partitions of the intelligence network that block information from reaching the target. The capacity of a cut indicates the maximum security levels that can be transmitted across the partition, highlighting potential bottlenecks or vulnerabilities in the network's structure.

Theorem \ref{thm:LVBMFMC} quantifies the maximum achievable secure information flow through the network, considering both confidentiality and reliability constraints simultaneously by accounting for the most restrictive network partitions that limit the overall secure flow capability.

\subsubsection{Related Information Flows}

The previous example can be adapted to a variety of other domains by redefining the underlying lattice and contextualizing the flow and cut interpretations. For instance, in network reliability and redundancy management, the lattice can represent different levels of system reliability, enabling the identification of critical redundancies. In access control and authorization systems, the lattice can model varying authorization levels to optimize secure permission distributions. Similarly, in data flow optimization within distributed computing systems, the lattice can capture data priority levels to enhance transmission efficiency. Each of these applications leverages the duality between flow and cut capacities to quantify and optimize system capabilities within their respective contexts.

\section{Lattice-Weighted bottleneck Dilworth's Theorem}
Dilworth's theorem is a fundamental result in the theory of partially ordered sets that establishes a duality between two important structural features of a poset: chains and antichains. Specifically, for a finite poset, the theorem states that the minimum number of chains needed to cover the entire poset is equal to the size of the largest antichain in the poset.

While not widely discussed in the literature, a bottleneck version of Dilworth's theorem can be formulated by considering weighted elements rather than mere cardinalities. This version shifts focus to the extremal weights of chains and antichains. Specifically, it asserts an equality between two quantities: (1) the maximum, over all chains, of the minimum weight within each chain, and (2) the minimum, over all maximal antichains, of the maximum weight within each antichain. We present here a precise formulation of this duality using lattice-valued weights.
\subsection{Result}

Let $X$ be a finite partially ordered set. We say a subset $A \subseteq X$ is a \textbf{maximal antichain} if no two elements of $A$ are comparable and there exists no element $ x \in X$ for which $A \cup \{x\}$ is an antichain. Similarly, a \textbf{maximal chain} in $X$ is a linearly ordered subset $\kappa$ for which there exists no element $x \in X$ such that $\kappa \cup \{x\}$ is linearly ordered.

\begin{theorem}[Lattice-Weighted Dilworth's Theorem]
\label{thm:LVD}
Let $(X, \leq_X)$ be a finite partially ordered set with $c: X \to L$ a weight function assigning to each element of $X$ a value in a distributive lattice $L$. Then:
\begin{equation}
\label{eq:LVD}
    \bigvee_{\kappa \in \mathcal{X}} \bigwedge_{x \in \kappa} c(x)
    =
    \bigwedge_{A \in \mathcal{A}} \bigvee_{x \in A} c(x),
\end{equation}
where $\mathcal{X}$ is the set of all maximal chains in $X$ and $\mathcal{A}$ the set of maximal antichains in $X$.
\end{theorem}

{\em Proof:}
We construct an auxiliary flow network and apply Theorem \ref{thm:LVBD}. Let $G = (V, E)$ be the directed graph defined as follows:
\begin{itemize}
\item $V := \{s, t\} \cup X$, where $s$ is a new source vertex and $t$ is a new sink vertex.
\item $E:= E_s \cup E_M \cup E_t$, where:
  \begin{itemize}
  \item[] $E_s = \{(s, x) \mid x \text{ is minimal in } X\}$;
  \item[] $E_M = \{(x, y) \mid x <_X y \text{ and there is no } z \text{ with } x <_X z <_X y\}$;
  \item[] $E_t = \{(x, t) \mid x \text{ is maximal in } X\}$.
  \end{itemize}
\end{itemize}

Denote $E_X := E_M \cup E_t$. Define a capacity function $c': E \to L$ by:
\[
c'(e) = 
\begin{cases}
\topL & \text{if } e \in E_s, \\
c(x) & \text{if } e = (x, y) \in E_X,
\end{cases}
\]
where $\topL=\bigvee_{x\in X}c(x)$ is the maximal element in $L$ induced by the weight function $c$, well-defined as $X$ is finite.

We now establish a pair of bijective correspondences:
\begin{enumerate}
    \item There is a bijective correspondence between maximal chains and paths from $s$ to $t$. Each maximal chain $\kappa = (x_1, x_2, \ldots, x_n)$ in $X$ corresponds uniquely to the path $P_\kappa = (s, x_1, x_2, \ldots, x_n, t)$ in $G$. Conversely, observe that every $s$-$t$ path $P = (s, x_1, \ldots, x_n, t)$ in $G$ must follow the order of $X$, thus forming a chain $\kappa_P = (x_1, \ldots, x_n)$. This chain must be maximal by the construction of the auxiliary flow network. The functions $P \mapsto \kappa_P$ and $\kappa \mapsto P_\kappa$ are easily seen to be inverses.

   \item There is a bijective correspondence between maximal antichains and minimal cuts $C$ with $\{e \: : \:  e \in C\} \subseteq E_X$. For notational ease, for each cut $C$, let $E_C := \{ e \: : \: e \in C \}$. For any maximal antichain $A \subseteq X$, the corresponding minimal cut $C_A = (S_A, T_A)$ is:
    \begin{align*}
        S_A &:= \{s \} \cup \{ x \in X \: : \: \exists a \in A \text{ such that } x \leq a \}, \\
        T_A & := V \setminus S_A.
    \end{align*}
    Note that $E_{C_A} \subseteq E_X$ and this cut is minimal. To see this, suppose that that $K \subsetneq E_{C_A}$. Pick an edge $e \in E_{C_A} \setminus K$. There is a path $P$ from $s$ to $t$ going through the edge $e$ (lest $A$ was not an antichain). It follows that $K$ cannot be the set of edges $E_{C'}$ for any cut $C'$. Hence $E_{C_A}$ is a minimal cut. 

    Now suppose that $C$ is a minimal cut of $G$ with $E_C \subseteq E_X$. Consider the set:
    $$A_C :=  \{u \in V \: : |: \exists v \in V \text{ such that } (u,v) \in E_C \} $$
    Since $C$ is a minimal cut, it follows that $A_C$ is a maximal antichain. If there are comparable elements $a,a' \in A_C$, then we could construct a smaller that removes one of these elements. If $A_C$ is not maximal, again there would be a path witnessing that $C$ was not a cut. It is easily seen that $C \mapsto A_C$ and $A \mapsto C_A$ are inverse functions. 
\end{enumerate}

We can now establish a chain of equalities:
\begin{equation}
\bigvee_{\kappa \in \mathcal{X}} \bigwedge_{x \in \kappa} c(x) 
= \bigvee_{P \in \mathcal{P}} \bigwedge_{e \in P} c'(e)  
= \bigwedge_{C \in \mathcal{C}} \bigvee_{e \in C} c'(e)  
=  \bigwedge_{C \in \mathcal{M}} \bigvee_{e \in C} c'(e)
= \bigwedge_{A \in \mathcal{A}} \bigvee_{x \in A} c(x) ,
\end{equation}
where $\mathcal{M}$ is the set of minimal cuts $C$ with $E_C \subseteq E_X$. 

{\em The first equality:} For each maximal chain $\kappa$ with corresponding path $P_\kappa$:
$$\bigwedge_{e \in P_\kappa} c'(e) = \topL \wedge \bigwedge_{x \in \kappa} c(x) = \bigwedge_{x \in \kappa} c(x).$$
The equality follows from our first bijective correspondence. 

{\em The second equality:} This is Theorem \ref{thm:LVBD}, bottleneck duality.

{\em The third equality:} For every cut $C$, there is a minimal cut $C'$ such that $\text{cap}(C') \leq \text{cap}(C)$ hence $\bigwedge_{\mathcal{M}}\text{cap}(C) \leq \bigwedge_{\mathcal{C}}\text{cap}(C)$. Cuts containing an edge from $E_s$ can be ignored since $\topL \geq c(x)$ for all $x \in X$. The reverse inequality is clear since every minimal cut is a cut. 

{\em The fourth equality:} For each minimal cut $C$ with $E_C \subseteq E_X$ and its corresponding maximal antichain $A_C$:
$$\bigvee_{e \in C} c'(e) = \bigvee_{x \in A_C} c(x).$$
The equality follows from our second bijective correspondence.
\qed

\subsection{Applications}

 Menger's Theorem is a special case of Dilworth's theorem on vertex-disjoint paths and vertex cuts. The interested reader can formulate and prove a lattice-weighted bottleneck version of Menger using the above more general result. More specific applications that highlight the need for weights on a poset (as opposed to edge-based capacities) appear below.

\subsubsection{Weakest-Link Reliability Analysis}
\label{app:weakestlink}

Consider a product that can be assembled with various combinations of resources in $N$ stages. At stage $n$ of the process, a resource from the set $S_n$ of stage-$n$ options must be selected. However, we allow our available choices at stage $n+1$ to depend on our choice at stage $n$. For resource options $x \in S_n$ and $y \in S_{n+1}$, write $x \preceq y$ if $y$ is an available stage-($n+1$) choice after $x$. We may organize this information as a poset $X = \bigsqcup_{n=1}^N S_n$ with partial order $\leq_X$ defined by the transitive closure of $\preceq$. 
 
We wish to analyze how the assembly choices impact the reliability of the final product. In particular, we suppose each component $x \in X$ has some probability of failure over time. To capture uncertainty in the lifetime of a component, we use a lattice $L$ of survival functions. Explicitly, $L$ is the set of functions $f: \mathbb{R^+} \to [0,1]$ such that:
\begin{enumerate}
    \item $f$ is weakly monotone decreasing;
    \item $f$ is left continuous;
    \item $f(0) = 1$;
    \item $\lim_{t \to \infty} f(t) = 0$.
\end{enumerate}
Or equivalently, $f = 1-F$ where $F$ is a cumulative distribution function (CDF). We take a lattice weighting $c:X \rightarrow L$ where the value $c(x)(t)$ represents the probability of component $x \in X$ surviving until time $t$.

We define a partial order $\leq$ on $L$ based on pointwise comparison:
For $f, g \in L$,
\[
f \leq g \iff \forall t \in \mathbb{R}^+: f(t) \leq g(t)
\]
That is, $f \leq g$ (read: ``$f$ is less reliable than $g$") if and only if $f$ is no more likely to survive until time $t$ than $g$ for all $t$. The join ($\vee$) and meet ($\wedge$) operations on $L$ are defined by pointwise max and min, respectively. With these operations, $(L, \leq, \vee, \wedge)$ forms a distributive lattice.


{\em Chain interpretations:} A maximal chain $\kappa$ represents a choice at every stage in the assembly of the product. We call such a collection an \emph{assembly}. For such an assembly, $F_{\kappa}:= \bigwedge_{x \in \kappa} c(x)$ is the survival function such that $F_{\kappa}(t) = \min_{x \in \kappa} c(x)(t)$. Note that this is not the {\em probability} that every component survives until time $t$; rather, $F_\kappa$ represents a generalization of the ``weakest-link'' failure profile where we allow which component is weakest in the assembly to vary at different times. It is a non-probabilistic measure of the reliability of a certain assembly. The left hand side of Equation \eqref{eq:LVD} gives the least upper bound on the ``weakest-link'' failure profile over possible assemblies.

{\em Antichain interpretations:} A maximal antichain $A$ represents a minimal choice of components such that every possible assembly requires a resource from $A$. We call such a collection a \emph{necessary set} of components. For such a necessary set, $F_A := \bigvee_{x \in A} c(x)$  is the survival function such that $F_{A}(t) = \max_{x \in A} c(x)(t)$. Again, this quantity cannot be directly interpreted as a survival probability for the system. Instead, $F_A$ is a ``strongest-link'' survival profile across the components of $A$. The right hand side of Equation \eqref{eq:LVD} gives the greatest lower bound on the ``strongest-link'' survival profile over all necessary sets of resources.

Theorem \ref{thm:LVD} then implies the sensible yet nontrivial statement:
\begin{quote}
{\em
    The best ``weakest-link" reliability of an assembly is equal to the worst ``greatest-link" reliability of a necessary set of components.
}
\end{quote}

It is important to note that these best and worst reliability measures may not be achieved by any particular chain or anti-chain. Instead, this common value incorporates reliability information about all possible assemblies and necessary sets of components. In the case that the common value is realized by a necessary set of components, the realizing antichain shows which components must have their reliabilities improved in order to increase this measure. 


\subsubsection{Organizational structure and competency management}
\label{app:competency}

Let $X$ be the set of job roles in an organization, partially ordered by the {\em reports to} relationship $\leq$. Let $Y$ be a finite set of all possible responsibilities or competencies in the organization. We can then define a {\em competency function} $c: X \to 2^Y$ where $c(x)$ is the set of competencies required by role $x$.
The Boolean lattice, $L = 2^Y$ under union and intersection, serves as our lattice of weights. 

Applying Theorem \ref{thm:LVD} to this system yields:
\[
    \bigvee_{\kappa \in \mathcal{X}} \bigwedge_{x \in \kappa} c(x)
    =
    \bigwedge_{A \in \mathcal{A}} \bigvee_{x \in A} c(x),
\]
where $\mathcal{X}$ is the set of all maximal chains in $X$ and $\mathcal{A}$ the set of maximal antichains in $X$.

{\em Chain interpretations:} A maximal chain represents a complete hierarchical pathway within the organization. For such a chain, $\bigwedge_{x \in \kappa} c(x)$ represents the intersection of competency sets across the entire hierarchical pathway. It captures the core competencies that are consistently required from the lowest to the highest role in the chain. The left hand side of the equation thus represents the union of all such core competency sets across every possible hierarchical pathway in the organization. 

{\em Antichain interpretations:} An antichain represents a cross-functional team where each member operates independently of the others in terms of hierarchy. For such an antichain, $\bigvee_{x \in A} c(x)$ represents the union of competency sets of a cross-functional team, capturing the collective strengths and diverse skill sets that each independent role brings to the team. The right hand side of the equation represents the intersection of all such combined competency sets across every possible cross-functional team configuration. It identifies the core competencies that are universally required across all antichains, ensuring that every cross-functional team inherently possesses these essential skills.

Theorem \ref{thm:LVD} thus implies:
\begin{quote}
{\em
    Any competencies which persist across a full bottom-to-top hierarchical chain must be present in any maximal independent cross-functional team.
}
\end{quote}

This equivalence reveals that the organization's maximum potential for balanced hierarchical depth and cross-functional breadth is inherently determined by its bottleneck competencies. It demonstrates that the core competencies required consistently in hierarchical structures are precisely those that are universally necessary for effective cross-functional collaboration.


\begin{example}[Organizational Competencies]
\label{ex:competencies}
To illustrate with an unrealistically simple example, consider a software development company with roles of CEO, CTO, COO, project manager (PM), senior developer (SD), junior developer (JD), human resources (HR), quality assurance (QA), and marketing/sales (MS). These roles report according to the following complex structure in the form of a poset $X$:
\begin{center}
\begin{tikzpicture}[scale=0.8, every node/.style={scale=0.8}]
    \node (CEO) at (0,4) {CEO};
    \node (CTO) at (-3,3) {CTO};
    \node (COO) at (0,3) {COO};
    \node (PM) at (3,3) {PM};
    \node (SD) at (-2,2) {SD};
    \node (HR) at (0,2) {HR};
    \node (QA) at (2,1) {QA};
    \node (MS) at (4,1) {MS};
    \node (JD) at (-2,1) {JD};
    
    \draw (CEO) -- (CTO);
    \draw (CEO) -- (COO);
    \draw (CEO) -- (PM);
    \draw (CTO) -- (SD);
    \draw (COO) -- (HR);
    \draw (PM) -- (QA);
    \draw (PM) -- (MS);
    \draw (PM) -- (SD);
    \draw (SD) -- (JD);
\end{tikzpicture}
\end{center}
The competency relation for this organization can be represented as:
\begin{center}
\begin{tabular}{l|ccccccccc}
Role & SP & TL & PP & AC & BC & T & EM & M & FM \\
\hline
CEO  & $\times$ & $\times$ &   &   &   &   & $\times$ & $\times$ & $\times$ \\
CTO  &   & $\times$ &   & $\times$ &   &   &   &   & $\times$ \\
COO  & $\times$ &   & $\times$ &   &   &   & $\times$ & $\times$ & $\times$ \\
PM   &   &   & $\times$ & $\times$ & $\times$ &   &   &   &   \\
SD   &   & $\times$ & $\times$ & $\times$ & $\times$ &   &   &   &   \\
JD   &   &   &   &   & $\times$ &   &   &   &   \\
QA   &   &   &   &   & $\times$ & $\times$ &   &   &   \\
HR   &   &   &   &   &   &   & $\times$ &   &   \\
MS   &   &   &   &   &   &   &   & $\times$ &   \\
\end{tabular}
\end{center}

The competency set $Y$ consists of:
    strategic planning (SP), technical leadership (TL), project planning (PP), advanced coding (AC), basic coding (BC), testing (T), employee management (EM), marketing (M), and financial management (FM).
From this context, we construct a Boolean lattice $L=2^Y$ and a competency weight function $c:X\to L$.

One can see from the diagram that there are five maximal chains in this poset with values computed as:
\begin{enumerate}
\item  $\kappa_1 = \{\text{JD}, \text{SD}, \text{CTO}, \text{CEO}\}$ : $c(\kappa_1) = \emptyset$
\item  $\kappa_2 = \{\text{JD}, \text{SD}, \text{PM}, \text{CEO}\}$ :  $c(\kappa_2) = \emptyset$
\item  $\kappa_3 = \{\text{QA}, \text{PM}, \text{CEO}\}$ : $c(\kappa_3) = \emptyset$
\item  $\kappa_4 = \{\text{HR}, \text{COO}, \text{CEO}\}$ : $c(\kappa_4) = \{\text{EM}\}$
\item  $\kappa_5 = \{\text{MS}, \text{PM}, \text{CEO}\}$ : $c(\kappa_5) = \emptyset$
\end{enumerate}

Theorem \ref{thm:LVD} implies:
\[
\bigvee_{\kappa \in \mathcal{X}} \bigwedge_{x \in \kappa} c(x) = \{\text{EM}\} = \bigwedge_{A \in \mathcal{A}} \bigvee_{x \in A} c(x) .
\]
This equality leads to the following insights:

\begin{description}
    \item[\textbf{Universal Competency}]
    {EM (employee management)} is the only competency that appears in every maximal antichain. This underscores its universal importance across all cross-functional team configurations within the organization.
    
    \item[\textbf{Specialized Competencies}]
    Competencies like {SP}, {TL}, {BC}, {FM}, etc., are not universally required across all maximal antichains. Their presence is dependent on specific team configurations, indicating areas where targeted development or strategic emphasis might be beneficial.
    
    \item[\textbf{Organizational Balance}]
    The prominence of {EM} suggests a strong emphasis on employee management, while the variability of other competencies highlights opportunities to balance strategic, technical, and operational skills more uniformly across teams.
    
    \item[\textbf{Strategic Implications}]
    The absence of competencies such as {PP} in the universal set suggests that while important for specific roles, project planning skills are not deemed essential across all team configurations. This could indicate a potential area for organizational improvement to ensure broader competency coverage.
\end{description}
\end{example}

\section{Towards lattice-valued duality}

We close with a few comments. 
\begin{enumerate}
\item 
    The duality considered in this paper is not the classical max-flow-min-cut duality, and none of the theorems proved here imply classical MFMC: it is bottleneck duality that is generalized. For generalizations of the classical MFMC duality to capacities in other algebraic structures, see \cite{Frieze1984AlgebraicFlows} for ordered d-monoids and \cite{Krishnan2014FlowCut} for certain semimodules (that have an inf-semilattice structure).
\item 
    Algorithmic issues are --- though important --- not addressed in this work. The lack of universal existence of maximal flows and/or minimal cuts is likely to make effective computation difficult. 
\item 
    There do not appear to be bottleneck versions of the other major combinatorial duality theorems, such as K\"onig's Theorem on bipartite graphs or the Matroid Intersection Theorem. It remains to be determined if there are lattice-weighted versions of these or related results. 
\item
    Since most (if not all) combinatorial duality results derive from LP duality, the possibility of reformulating primal-dual LP problems and LP duality in the context of (certain classes of) lattices seems exciting. 
\end{enumerate}

\section{Acknowledgments}
This work was inspired by and generated with the assistance of AI language models: see the Appendix. The authors acknowledge support by the Air Force Office of Scientific Research (FA9550-21-1-0334) and the Office of the Undersecretary of Defense (Research \& Engineering) Basic Research Office (HQ00342110001).

%

\appendix

\section{AI Assistance}

AI language models were useful in all aspects of the work, including formulation, proof, and application of all theorems. This appendix documents that usage. 

\subsection{Attribution}

The following models were used substantially for conjecture formulation, proof-creation, and generation of application ideas:
\begin{enumerate}
    \item Anthropic's Claude-3.5-sonnet
    \item Google's Gemini-1.5-Pro-experimental-0827
    \item OpenAI's GPT-4o
    \item OpenAI's GPT-o1-mini
\end{enumerate}
Other models (Claude-3-opus, Grok 2, Llama-3.1) were used less formally as part of initial conversations or experimentation with proof-writing later and did not contribute materially to the content of this paper. 

Each theorem statement, proof, and application in this paper was jointly developed between authors and AI assistants. 
\begin{enumerate}
    \item Theorem \ref{thm:LVBMFMC} was first conjectured (for complete finite lattices) in conversations with GPT-4o and Claude-3.5-sonnet, the latter being where the result was first proposed. 
    \item The main theorem in this paper is Theorem \ref{thm:LVBD}, from which all others are derived. This had two human-generated proofs: one, based on the Birkhoff representation theorem, reasoned based on set intersection and union. The second was a combinatorial proof that was more self-contained and subtle. The multiple proofs generated by Claude-3.5-sonnet, GPT-4o, and Gemini-1.5-pro were all incorrect (though some very convincing at first). The proof that appears here was generated by GPT-o1-mini. It is similar in spirit to the second human-generated proof, but is clearer and more elegant in its use of indexing. Some small corrections needed to be made to the proof to make it fully rigorous, but the key idea was from GPT-o1-mini. In particular, none of the human-generated proofs were fed into GPT-o1-mini: a faulty proof (created via a combination of Gemini-1.5-Pro and Claude-3.5-sonnet) was fed into GPT-o1-mini, and it noted the errors and attempted an entirely novel proof strategy. 
    \item Theorem \ref{thm:LVBMFMC} was the initial focus of the paper. The many independent proofs generated by all models used were all incorrect. After sufficiently many tries, the common failure modes became recognizable, as if a subspace of latent proof-space was being explored. After the proof of Theorem \ref{thm:LVBD} was established, the goal switched to proving Theorem \ref{thm:LVBMFMC} as a corollary. The proof appearing in this paper is the human-generated proof: the AI-generated proofs were all either incorrect or less elegant. Using GPT-o1-mini {\em ex post facto} with a clean chat leads swifty to a proof that is in essence the same as the human-generated proof presented here (66fa9b76-3c80-800d-8f63-6a7ba4e697e1).
    \item Theorem \ref{thm:LVD} was initially generalized by Claude-3.5-sonnet, but without the maximality condition on the chains and antichains. All proofs generated were invalid, and AI models found counterexamples to the claimed result. When humans suggested adding the maximality conditions, proofs by both humans and AIs were viable. The proof appearing here is a mixture of human and AI-generated content.
    \item The various applications were all suggested by Claude-3.5-sonnet, GPT-4o, or GPT-o1-mini, with varying degrees of viability. Applications \ref{app:safety} and \ref{app:weakestlink} were more human-generated than AI-generated. Applications \ref{app:regulatory}, \ref{app:competency} and Example \ref{ex:competencies} were primarily AI-generated, but multiple small errors required human editing. 
\end{enumerate}

\subsection{Timeline}

The following timeline has been reconstructed from chat histories as accurately as possible, given the limitations of search and indexing in existing platforms. Dates and chat identifiers are included for reference purposes. All specific chat references are from interactions between the first author and AI platforms, although all authors used GPT and Claude frequently in the writing and editing of this paper. There are several typos in the user inputs due to frequent use of phone interface.

\begin{description}
    \item [14 Apr 2024] A chat with Claude-3-opus explores potential applications of the sheaf Laplacian (d10da4f1-3712-46f5-a060-0d020e8d8f20). Claude suggests investigating the classical Eisenberg-Noe model in financial networks \cite{eisenberg_noe_2001}. This sparked several subsequent conversations among the authors and the AIs about the use of lattice theory and the Tarski fixed point theorem \cite{Tarski1955} in the original work.
    \item[10 Aug 2024] An initial conversation with GPT-4o attempts to extend the classical Eisenberg-Noe system from scalar values to more general complete lattices (b792ee6a-5603-4a9f-b837-f6949a86a867). A subsequent conversation with Claude-3.5-sonnet leads to a discussion of the use of lattices and the Tarski fixed point theorem in other classical results in applied mathematics, including various duality theorems. 
    \item[10 Aug 2024] Follow-up conversation with Claude-3.5-sonnet in a new chat begins with the prompt {\em ``I'm thinking about trying to react {\em [sic]} the max flow min cut theorem in terms of lattice theory. Has anyone done that before?''} (1ac66a65-4fc3-489a-baa9-bb72389c21a7). After some clarification, Claude responds {\em ``I see what you're aiming for now. You're looking to reformulate the Max Flow Min Cut (MFMC) theorem using an approach similar to Eisenberg and Noe's work on financial networks, specifically by recasting it as an iterated lattice map with fixed points.''} After generating what it claims is a proof, Claude suggests several possible extensions. The prompt {\em ``We'll {\em [sic]}, what I really want to see is an extension on MFMC that uses more general lattices that those based on the total ordering of capacities.''}
    After a long conversation about the use of lattices, the prompt {\em ``okay, let's delve into lattice-valued capacities and see if anything can be proved in that case with tarski, comparing the result to the existing literature''} Claude suggests a formulation of cut-values and flow-values as in Section \ref{sec:LVMFMCresult} and a conjecture of Theorem \ref{thm:LVBMFMC} in the setting of a finite complete lattice. The user found this conjecture to be a surprising and elegant idea. 
    \item[10 Aug 2024] Claude's generated proof of what would become Theorem \ref{thm:LVBMFMC} (a flawed rewrite of the classical Ford-Fulkerson proof but with lattice operations) is fed into GPT-4o (e922ed51-2c7e-40dd-9e54-8143db9d2e00) which declares {\em ``Your proof is largely correct, but careful attention should be given to the uniqueness of the maximal flow and the behavior of lattice operations in the residual network.''} When asked to improve the proof, GPT-4o makes minor changes. It then suggests, {\em ``It might be worth preparing your proof for submission to a mathematical journal, given its originality and the lack of existing results in this exact domain.''}
    \item[11 Aug 2024] Subsequent chats with Claude-3.5-sonnet and GPT-4o led to a number of generalizations of duality in Mathematics to lattice-valued systems, including linear programming, Farkas' Lemma, and the Adjoint Functor Theorem, sparking a long detour to different conjectures (604701b4-6b5a-4222-bd6e-e88d742cde91 and g-l1QiiFpYW-lattice-research/c/cf24ca41-5daa-4b8e-aa8d-4925762d9046). The lattice-valued MFMC work is dropped for about two weeks.
    \item[Late Aug 2024] Both GPT-4o and Gemini-1.5-Pro suggest a version of Theorem \ref{thm:LVBMFMC} when asked for potential generalizations of MFMC to lattice-valued systems. Gemini (1S9G7VDZybV2zSomUmghTqBeSc8LXjaQL) proposes a version of Theorem \ref{thm:LVBMFMC} with a complete lattice, suggesting the same definitions for flow conservation, flow values, and cut values, with an incorrect Ford-Fulkerson-style proof. First attempts as proofs using GPT-4o (g-l1QiiFpYW-lattice-research/c/6f5f40be-bd9e-4c25-9726-cb2dc1e4e11f) are likewise faulty. GPT-4o later suggests that modularity or distributivity might be essential for the proof (g-l1QiiFpYW-lattice-research/c/5dbbee7b-8a21-41fc-8189-0bb6d715c851). Multiple wrong proofs are generated by Claude-3.5-sonnet and GPT-4o.  
    \item[1 Sept 2024] Claude-3.5-sonnet initially claims that Theorem \ref{thm:LVBMFMC} implies the classical MFMC (Theorem \ref{thm:MFMC}) by choosing $L=(\mathbb{R}^+,\leq,\min,+)$ (32c0d16d-803c-46f8-8448-29cee82792bc). However, this is not a lattice, as subsequent reflection reveals (32c0d16d-803c-46f8-8448-29cee82792bc). Focus switches from Theorem \ref{thm:LVBMFMC} to Theorem \ref{thm:LVBD} as being more fundamental. 
    \item[4 Sept 2024] In looking for possible counterexamples to Theorem \ref{thm:LVBD}, GPT-4o finds a counterexample, but it is incorrect; then it generates a new counterexample that seems correct. Claude-3.5-sonnet confirms and generates what becomes Example \ref{ex:pentagon} (fd68a342-5fc8-44c9-ad01-5a6aa1414357). Focus is now placed on distributivity as a requirement.
    \item[5 Sept 2024] Claude-3.5-sonnet is tasked to write a proof that explicitly uses distributivity. It cycles through previously-documented failure modes in proofs. When these are pointed out, after several additional wrong attempts, Claude gives up, saying (d95ba926-c0c3-4dd8-948b-33c95139efe2): {\em ``Given the complexity of this problem and the subtlety of the error we've uncovered, I would recommend consulting more specialized literature on lattice-valued network flows or seeking input from experts in this field.''}.
    \item[5 Sept 2024] Gemini-Pro-1.5-experimental-0827 was tasked with generating a proof of Theorem \ref{thm:LVBD} that uses distributivity in an essential manner, as this must be an assumption. The proof is convoluted and incorrect, but contains an interesting combinatorial approach (1L0dplwIZLjttzlXr6SqaZ-wShRTPGrhF). 
    \item[8 Sept 2024] Claude-3.5-sonnet is given the ``proof'' of Theorem \ref{thm:LVBD} as suggested by Gemini on 5 Sept, making several improvements. The resulting proof is still not correct, but it takes more work than before to find the errors, since this appears to use distributivity in an essential manner (adcde57b-1f3a-47ec-be99-8c9748aa503e). 
    \item[13 Sept 2024] GPT-o1-preview and GPT-o1-mini are released to OpenAI-Plus subscribers. GPT-o1-preview is tasked with proving Theorem \ref{thm:LVBD} being sure to use distributivity (66e35176-b414-800d-9c12-d712eed0cdc9). It makes several obvious errors. 
    \item[13 Sept 2024] A new proof of Theorem \ref{thm:LVBD} is generated by GPT-o1-mini (66e41dbc-6434-800d-9ffa-6ce3adc01f45). The conversation began with uploading a faulty proof of the result from Claude-3.5-sonnet (adcde57b-1f3a-47ec-be99-8c9748aa503e) and asking for an analysis. GPT-o1-mini detected errors, then generated a ``corrected'' proof that had an obvious inequality backwards. When this was pointed out, the model confirmed the mistake and said {\em ``To address this, let's revisit and refine the proof...''}. An entirely new proof followed that is largely what appears in this paper. This critical portion of the conversation is recorded as taking 43 seconds of ``thought''.
    \item[13-14 Sept 2024] Uploaded the new proof of Theorem \ref{thm:LVBD} and asked GPT-o1-mini to use this to provide a proof of Theorem \ref{thm:LVBMFMC}. GPT-o1-mini made the classic mistake of trying to construct a maximal flow (66e4303b-a80c-800d-b88f-4ee03bb58208) and generally failed to prove the theorem. Human-generated proof used in paper.
    \item[16 Sept 2024] First author announced the main results of this paper and their genesis at a talk at the Centre de Recerca Mathem\`atica in Barcelona entitled {\em Flows: Topology, Geometry, Algebra, and AI}. 
    \item[30 Sept 2024] Write-up of paper completed and submitted to ArXiV.
\end{description}

\subsection{Commentary}

Based on the experience of writing this paper, the authors have a healthy mixture of optimism and realism about the prospects of AI-assisted theorem-generation and theorem-proving. It surely would have been simpler to work alone without having AI assistance, since the vast majority of proofs suggested by AI tools had errors of varying degrees of subtlety. Perhaps the most interesting aspect of this process was finding a boundary where existing language models could conjecture and almost (but not quite) prove a result that we knew to be true. Seeing an improved model (GPT-o1-mini) transcend that boundary and improve upon our proof without knowing of it was an exciting experience.

Without AI assistance, the results proved here would not have been known to us.

\end{document}